\documentclass{amsart}
\usepackage[alphabetic]{amsrefs}
\usepackage{amssymb,amscd}
\usepackage[all]{xy}

\numberwithin{equation}{subsection}

\swapnumbers
\theoremstyle{plain}
\newtheorem{thm}[subsection]{Theorem}
\newtheorem{prop}[subsection]{Proposition}
\newtheorem{cor}[subsection]{Corollary}
\newtheorem*{thm*}{Theorem}

\theoremstyle{definition}

\theoremstyle{remark}

\newtheorem{rems}[subsection]{Remarks}

\newcommand{\Curve}{\mathcal{C}}
\renewcommand{\O}{\mathcal{O}}

\newcommand{\F}{\mathbb{F}}
\newcommand{\Fp}{{\mathbb{F}_p}}

\newcommand{\Fq}{{\mathbb{F}_q}}
\newcommand{\Fqtimes}{{\mathbb{F}_q^\times}}
\newcommand{\Z}{\mathbb{Z}}
\newcommand{\Q}{\mathbb{Q}}

\renewcommand{\P}{\mathbb{P}}

\renewcommand{\d}{\mathfrak{d}}

\newcommand{\n}{\mathfrak{n}}
\newcommand{\p}{\mathfrak{p}}

\newcommand{\sha}{{\hbox to 10pt{\rlap{\hskip2.8pt\vrule
height6pt\hskip1.6pt\vrule height6pt\hskip1.6pt
\vrule height6pt}\hskip1pt\vrule height0.8pt width 8pt\hskip1pt}}}

\newcommand{\into}{\hookrightarrow}

\newcommand{\tensor}{\otimes}
\newcommand{\compose}{\circ}
\newcommand{\nodiv}{\not|}
\def\nodiv{\mathrel{\mathchoice{\not|}{\not|}{\kern-.2em\not\kern.2em|}{\kern-.2em\not\kern.2em|}}}

\DeclareMathOperator{\im}{Im}
\DeclareMathOperator{\tr}{Tr}
\DeclareMathOperator{\N}{N}

\DeclareMathOperator{\ord}{ord}
\DeclareMathOperator{\rk}{Rank}

\DeclareMathOperator{\Hom}{Hom}

\DeclareMathOperator{\gal}{Gal}
\DeclareMathOperator{\spec}{Spec}
\DeclareMathOperator{\en}{End}
\DeclareMathOperator{\mor}{Mor}

\renewcommand{\vector}[1]{\mathbf{#1}}
\renewcommand{\v}{\vector{v}}

\newcommand{\Fpbar}{{\overline{\mathbb{F}}_p}}
\newcommand{\Fqbar}{{\overline{\mathbb{F}}_q}}
\newcommand{\Fpf}{{\mathbb{F}}_{p^f}}
\newcommand{\Fr}{{\mathbb{F}}_{r}}
\newcommand{\Qbar}{{\overline{\mathbb{Q}}}}

\newcommand{\Ql}{{\mathbb{Q}_\ell}}
\newcommand{\Zl}{{\mathbb{Z}_\ell}}

\DeclareMathOperator{\sel}{Sel}

\newcommand{\Zhatp}{{\hat\Z^{(p)}}}

\renewcommand{\a}{{\mathbf{a}}}

\theoremstyle{plain}
\newtheorem{subsublemma}[subsubsection]{Lemma}

\begin{document}
\title[Small ranks over function fields]{Jacobi sums, Fermat Jacobians, 
\\and ranks of abelian varieties 
\\over towers of function fields} 
\author{Douglas Ulmer}
\address{Department of Mathematics \\ University of Arizona \\ Tucson,
AZ  85721}
\email{ulmer@math.arizona.edu}
\thanks{This paper is based upon work supported by the National
Science Foundation under Grant No. DMS 0400877}
\date{September 23, 2006}
\subjclass[2000]{Primary 14G05, 11G40; 
Secondary 11G05, 11G10, 11G30, 14G10, 14G25, 14K12, 14K15}
\maketitle

\section{Introduction}

\subsection{}
Given an abelian variety $A$ over a function field $K=k(\Curve)$ with
$\Curve$ an absolutely irreducible, smooth, proper curve over a field
$k$, it is natural to ask about the behavior of the Mordell-Weil group
of $A$ in the layers of a tower of fields over $K$.  The simplest
case, which is already very interesting, is when $A$ is an elliptic
curve, $K=k(t)$ is a rational function field, and one considers the
towers $k(t^{1/d})$ or $\overline k(t^{1/d})$ as $d$ varies through
powers of a prime or through all integers not divisible by the
characteristic of $k$.

When $k=\Q$ or more generally a number field, several authors (e.g.,
\cite{Shioda}, \cite{Stiller}, \cite{Fastenberg}, \cite{Silverman1},
\cite{Silverman2}, and \cite{Ellenberg}) have considered this question
and given bounds on the rank of $A$ over $\Q(t^{1/d})$ or
$\Qbar(t^{1/d})$.  In some interesting cases it can be shown that $A$
has rank bounded independently of $d$ in the tower $\Qbar(t^{1/d})$.
Of course no example is yet known of an elliptic curve over $\Q(t)$
with unbounded ranks in the tower $\Q(t^{1/d})$, nor of an elliptic
curve over $\Qbar(t)$ with non-constant $j$-invariant and unbounded
ranks in the tower $\Qbar(t^{1/d})$.

When $k$ is a finite field, examples of Shioda and the author show
that there are non-isotrivial elliptic curves over $\Fp(t)$ with
unbounded ranks in the towers $\Fpbar(t^{1/d})$ \cite{Shioda}*{Remark
  10} and $\Fp(t^{1/d})$ \cite{UlmerR1}*{1.5}.  More recently, the
author has shown \cite{UlmerR2} that high ranks over function fields
over finite fields are in some sense ubiquitous.  For example, for
every prime $p$ and every integer $g>0$ there are absolutely simple
abelian varieties of dimension $g$ over $\Fp(t)$ with unbounded ranks
in the tower $\Fp(t^{1/d})$, and given any non-isotrivial elliptic
curve $E$ over $\Fq(t)$, there exists a finite extension $\Fr(u)$ such
that $E$ has unbounded (analytic) ranks in the tower $\Fr(u^{1/d})$.

One obvious difference between number fields and finite fields which
might be relevant here is the complexity of their absolute Galois
groups: that of a finite field is pro-cyclic while that of a number
field is highly non-abelian.  Ellenberg uses this non-abelianess in a
serious way in his work on bounding ranks and, in a private
communication, he asked whether it might be the case that, say, a
non-isotrivial elliptic curve over $\Fq(t)$ always has unbounded rank
in the tower $\Fp(t^{1/d})$.

Our goal in this note, which is a companion to \cite{UlmerR2}, is to
give a number of examples of abelian varieties over function fields
$\Fq(t)$ which have bounded ranks in the towers $\Fqbar(t^{1/d})$ as
$d$ ranges through powers of a suitable prime or through all integers
not divisible by $p$, the characteristic of $\Fq$.  We also get some
information about ranks in towers $k(t^{1/d})$ for arbitrary fields
$k$.  Along the way we prove some new results on Fermat curves which
may be of independent interest.  The main results are
Theorems~\ref{thm:GaussJacobiBounds}, \ref{thm:Jd0}, \ref{thm:Jdp},
\ref{thm:isotrivial}, \ref{thm:nonisol}, and \ref{thm:nonisod}

\subsection{}
It is a pleasure to thank Jordan Ellenberg for his stimulating
questions about ranks of elliptic curves as well as Brian Conrey, Bill
McCallum, and Dinesh Thakur for their help.  Special thanks are due to
Bjorn Poonen for several incisive remarks and for pointing out that
some arguments originally given for elliptic curves apply more
generally to higher-dimensional abelian varieties.

\section{Jacobi sums}
 
\subsection{} 
Throughout the paper $p$ will be a rational prime number, $\Fp=\Z/p\Z$
will be the prime field of characteristic $p$, and $q=p^f$ will be a
power of $p$.  Fix an algebraic closure $\Qbar$ of $\Q$.  All number
fields considered will tacitly be assumed to be subfields of $\Qbar$.
We denote by $\mu_d$ the group of $d$-th roots of unity in $\Qbar$.

Let $\p$ be a prime of $\O_\Qbar$, the ring of integers of $\Qbar$,
over $p$.  The field $\O_\Qbar/\p$ is an algebraic closure of $\Fp$
which we denote by $\Fpbar$ and we write $\Fq$ for its subfield of
cardinality $q$.

Reduction modulo $\p$ induces an isomorphism between the group of all
roots of unity of order prime to $p$ in $\O_\Qbar$ and the
multiplicative group $(\O_\Qbar/\p)^\times=\Fpbar^\times$.  We let
$t:\Fpbar^\times\to\Qbar^\times$ denote the inverse of this
isomorphism.  We will use the same letter $t$ for the restriction to
any of the finite fields $\Fqtimes$.  Every character of $\Fqtimes$ is
a power of $t$.

\subsection{}
Fix a non-trivial additive character $\psi_p:\Fp\to\Qbar^\times$.  For
each $q$ we define an additive character $\psi_q$ as
$\psi_q=\psi_p\compose\tr_{\Fq/\Fp}$.
 
For each $q$ and each character $\chi$ of $\Fqtimes$, we define a
Gauss sum
$$G_q(\chi)=-\sum_{x\in\Fqtimes}\chi(x)\psi_q(x)\in\Q(\mu_{p(q-1)}).$$
It is well known that $G_q(\chi)=1$ if $\chi$ is the trivial character
and that $G_q(\chi)$ is an algebraic integer with absolute value
$q^{1/2}$ in every complex embedding if $\chi\neq1$.

For $d$ prime to $q$, $a\in\Z/d\Z$, and any $q\equiv1\pmod d$ we write
$G_q(a)$ for $G_q(t^{-a(q-1)/d})$ which lies in $\Q(\mu_{pd})$.  The
analysis leading to Stickelberger's theorem \cite{Washington}*{6.2}
shows that if $\wp$ is the prime of $\Q(\mu_{pd})$ under $\p$,
$q=p^f$, and $a\not\equiv0\pmod{d}$ then
$$\ord_{\wp}G_q(a)=(p-1)\sum_{j=0}^{f-1}
\left\langle\frac{p^ja}{d}\right\rangle$$
where $\langle x\rangle$ is the fractional part of $x$, i.e.,
$0\le\langle x\rangle<1$ and $x-\langle x\rangle\in\Z$.

\subsection{}\label{ss:JacobiSums}
Fix a positive integer $w$.  For each $q$ and each tuple of
non-trivial characters $\chi_0,\dots,\chi_{w+1}$ of $\Fqtimes$ such
that the product $\chi_0\cdots\chi_{w+1}$ is trivial, we define a
Jacobi sum
$$J_q(\chi_0,\dots,\chi_{w+1})=\frac1{q-1}
\sum_{\substack{x_0,\dots,x_{w+1}\in\Fpf^\times\\x_0+\cdots+x_{w+1}=0}}
\chi_0(x_0)\cdots\chi_{w+1}(x_{w+1})\in\Q(\mu_{q-1}).$$

It is well-known and elementary (see \cite{WeilNS}*{p.~501} for
example) that
$$J_q(\chi_0,\dots,\chi_{w+1})=\frac{(-1)^w}q\prod_{i=0}^{w+1}G_q(\chi_i).$$
In particular, the Jacobi sum is an algebraic integer with absolute
value $q^{w/2}$ in every complex embedding.

Let $A_{d,w}\subset(\Z/d\Z)^{w+2}$ be the set of tuples
$\a=(a_0,\dots,a_{w+1})$ such that $a_i\neq0$ for all $i$ and $\sum
a_i=0$.  If $\a\in A_{d,w}$ and $q\equiv1\pmod d$, we write $J_q(\a)$
for $J_q(t^{-a_0(q-1)/d},\dots,t^{-a_{w+1}(q-1)/d})$; clearly
$J_q(\a)\in\Q(\mu_d)$.  If $\wp$ is the prime of $\Q(\mu_{d})$ under
$\p$ and $q=p^f$, then
$$\ord_\wp J_q(\a)=\sum_{i=0}^{w+1}\sum_{j=0}^{f-1}
\left\langle\frac{p^ja_i}{d}\right\rangle-f.$$
 
We write $A'_{d,w}$ for those $\a\in A_{d,w}$ such that
$\gcd(d,a_0,\dots,a_{w+1})=1$.  Note that if $\a\in A_{d,w}$ and if
$e=\gcd(d,a_0,\dots,a_{w+1})$, $d'=d/e$ and
$\a'=(a_0/e,\dots,a_{w+1}/e)\in A'_{d',w}$ then for any $q\equiv1\pmod
d$ we have $J_q(\a)=J_q(\a')$.

Many of our results on ranks will be based on part (2) of the
following theorem about the distribution of Gauss and Jacobi sums.
Roughly speaking, it says that sums involving characters of large
order must either have large degree over $\Q$ or have valuation
bounded away from 0.

\begin{thm}\label{thm:GaussJacobiBounds}
\hfill\break
\begin{enumerate}
\item \vskip-12pt Fix a real number $\epsilon>0$ and a positive
  integer $n$.  There exists a constant $C_{\epsilon,n}$ depending
  only on $\epsilon$ and $n$ such that if $d>C_{\epsilon,n}$,
  $q=p^f\equiv1\pmod d$, $a\in(\Z/d\Z)^\times$, and the degree of
  $G_q(a)$ over $\Q(\mu_p)$ is $\le n$, then
$$\left|\frac{\ord_\wp G_q(a)}{(p-1)f}-\frac12\right|<\epsilon.$$ 
Here $\wp$ is the prime of $\Q(\mu_{pd})$ under $\p$.  Note that
  $\ord_\wp(q)=(p-1)f$.
\item Fix a positive integer $n$.  There exist constants $C_n$ and
  $\epsilon_n>0$ depending only on $n$ such that if $d>C_n$,
  $q=p^f\equiv1\pmod d$, $w\ge1$, $\a\in A'_{d,w}$, and the degree of
  $J_q(\a)$ over $\Q$ is $\le n$, then
$$\frac{\ord_\wp J_q(\a)}{f}>\epsilon_n.$$
Here $\wp$ is the prime of $\Q(\mu_d)$ under $\p$.  Note that
$\ord_\wp(q)=f$.
\end{enumerate}
\end{thm}
 
\begin{rems}
\hfill\break
\begin{enumerate}
\item \vskip-12pt The constants appearing in the theorem are {\it
    independent of\/} $p$ and effectively computable.
\item In part (2) of the theorem, we may replace ``the degree of
  $J_q(\a)$ over $\Q$ is $\le n$'' with ``the degree of the largest
  subfield of $\Q(J_q(\a))$ in which $p$ splits completely is $\le
  n$'' and similarly in part (1).  I do not know whether this has any
  applications to geometry.
\end{enumerate}
\end{rems}

The theorem is a consequence of Stickelberger's theorem and the
following simple estimate.
 
\begin{prop}\label{prop:BasicEstimate}
  Fix a real number $\epsilon>0$ and a positive integer $n$.  There
  exists a constant $C_{\epsilon,n}$ depending only on $\epsilon$ and
  $n$ such that if $d>C_{\epsilon,n}$ and $H\subset G=(\Z/d\Z)^\times$
  is a subgroup of index $\le n$, then for all $a\in G$,
$$\left|\frac{1}{|H|}\sum_{t\in H}\left\langle\frac{ta}{d}
\right\rangle-\frac12\right|< \epsilon.$$
\end{prop}
 
\begin{proof}
We have
\begin{align*}
 A:=\frac{1}{|H|}\sum_{t\in H}\left\langle\frac{ta}{d}\right\rangle
 &=\frac{1}{|H|}\sum_{\substack{s=1\\(s,d)=1}}^{d-1}\frac sd\frac1{[G:H]}
 \sum_{\chi\in\widehat{G/H}}\chi(sa^{-1})\\
 &=\frac12 +\frac1{d\phi(d)}\sum_{1\neq\chi\in\widehat{G/H}}\chi(a^{-1})
 \sum_{\substack{s=1\\(s,d)=1}}^{d-1}\chi(s)s
\end{align*}
where $|H|$ denotes the order of $H$, $\widehat{G/H}$ denotes the
group of characters of $G/H$ (which we view as characters of $G$
trivial on $H$), and $\phi(d)=|G|$ is Euler's function.  Partial
summation and the Polya-Vinogradov inequality \cite{Davenport}*{\S23}
show that there is an absolute constant $C$ such that the inner sum
above is $<Cd^{3/2}\log d$ and so the quantity $A$ to be estimated
satisfies
 $$\left|A-\frac12\right|\le\frac{Cnd^{1/2}\log d}{\phi(d)}.$$
 Well-known estimates for $\phi(d)$ \cite{HardyWright}*{Thm~327} say
 that for all $\delta>0$, $\phi(d)/d^{1-\delta}\to\infty$ as
 $d\to\infty$ so there is a constant $C_{\epsilon,n}$ depending only
 on $n$ and $\epsilon $ such that
 $$\frac{Cnd^{1/2}\log d}{\phi(d)}< \epsilon$$
whenever $d>C_{\epsilon,n}$.  This completes the proof of the proposition.
\end{proof}
 
\begin{cor}\label{cor:alld}
  Given $n$ there exists a constant $\delta_n>0$ depending only on $n$
  such that for any $d\ge2$, any $0\neq a\in\Z/d\Z$, and any subgroup
  $H\subset G=(\Z/d\Z)^\times$ of index $\le n$,
$$\frac{1}{|H|}\sum_{t\in H}\left\langle\frac{ta}{d}\right\rangle>\delta_n$$
 \end{cor}
 
\begin{proof}
  For $0\neq a\in(\Z/d\Z)$, set $e=\gcd(a,d)$, $d'=d/e$,
  $G'=(\Z/d'\Z)^\times$, $a'=a/e$, and $H'=\im(H\to G')$.  Then the
  index of $H'$ in $G'$ is $\le n$ and we have
 $$A:=\frac{1}{|H|}\sum_{t\in H}\left\langle\frac{ta}{d}\right\rangle=
 \frac{1}{|H'|}\sum_{t\in H'}\left\langle\frac{ta'}{d'}\right\rangle$$
and so we may assume that $\gcd(a,d)=1$, i.e., that $a\in G$.
 
Given $n$, let $C_{1/4,n}$ be the constant furnished by the
proposition for $n$ and $\epsilon =1/4$.  If $d>C_{1/4,n}$ then by the
proposition, $A>1/4$.  On the other hand, there are only finitely many
$d\le C_{1/4,n}$ and for each $d$, only finitely many subgroups
$H\subset(\Z/d\Z)$ of index $\le n$.  Since $A>0$ for each of these
finitely many possibilities, there is a $\delta_n>0$ such that
$A>\delta_n$ for all $d$ and $a$.
 \end{proof}

\subsection{Proof of Theorem~\ref{thm:GaussJacobiBounds} (1)}
Given $\epsilon$ and $n$, suppose that $d$, $q=p^f\equiv1\pmod d$, and
$a\in G=(\Z/d\Z)^\times\cong\gal(\Q(\mu_{pd})/\Q(\mu_p))$, are such
that $G_q(a)\in\Q(\mu_{pd})$ has degree $\le n$ over $\Q(\mu_p)$.  Let
$H\subset G$ be the subgroup of $G$ fixing $\Q(\mu_p,G_q(a))$, so that
$H$ has index $\le n$ in $G$.  If $\wp$ is the prime of $\Q(\mu_{pd})$
under $\p$, then for every $t\in H$, we have
$\ord_{\wp^t}(G_q(a))=\ord_{\wp}(G_q(a))$.  Therefore,
 \begin{align*}
 \frac{\ord_{\wp}G_q(a)}{(p-1)f}
 &=\frac1{|H|}\sum_{t\in H} \frac{\ord_{\wp^t}G_q(a)}{(p-1)f}\\
 &=\frac1{|H|f}\sum_{t\in H}\sum_{j=0}^{f-1}
\left\langle \frac{p^jta}{d}\right\rangle
 \end{align*}
 where the second equality comes from Stickelberger's theorem.  Let
 $P$ be the subgroup of $(\Z/d\Z)^\times$ generated by $p$ and $HP$
 the subgroup generated by $H$ and $P$.  The last displayed sum is
 then equal to
 $$\frac{1}{|HP|}\sum_{t\in
  HP}\left\langle\frac{ta}{d}\right\rangle.$$ 
Since $H$ has index $\le n$ in $G$, the same is true of $HP$ and so
Proposition~\ref{prop:BasicEstimate} shows that if $d>C_{\epsilon,n}$
then
$$\left|\frac{\ord_{\wp}G_q(a)}{(p-1)f}-\frac12\right|<\epsilon$$
as was to be shown.
\qed

\subsection{Proof of Theorem~\ref{thm:GaussJacobiBounds} (2)}
Given $\a\in A'_{d,w}$, set $d_i=d/\gcd(d,a_i)$.  The following lemma
tells us that if $d$ is large then at least two of the $d_i$ are also
large.
 
\begin{subsublemma} 
With notation as above, there exists an absolute constant $C$ such
that at least two of the $d_i$ are $\ge C\log d$.
\end{subsublemma}
 
\begin{proof}
If $\ell$ divides $d$ then from the definitions, there are at least
two $i$'s such that $\ell$ does not divide $a_i$.  Therefore the
largest prime power dividing $d$ also divides at least two of the
$d_i$.
 
To finish we note that Chebyschev's theorem
\cite{HardyWright}*{Thm.~7} implies that the the largest prime power
dividing $d$ is $\ge C'\log d$ for some absolute constant $C'$.
Indeed, let $M$ be a positive number, let $p_1,\dots,p_{\pi(M)}$ be
the primes less than $M$, and let $p_i^{e_i}$ be the largest power of
$p_i$ less than $M$.  If $N=\prod_{i=1}^{\pi(M)}p_i^{e_i}$ then
 $$\log N=\sum_{i=1}^{\pi(M)}e_i\log p_i\le\pi(M)\log M\le C'M$$ 
by Chebyschev.  This shows that if $N$ is a product of prime powers
less than $M$, then $N\le e^{C'M}$.  Therefore the largest prime power
dividing $N$ is at least $C\log N$ where $C=1/C'$.
\end{proof}

Now fix $n$ and consider those $q\equiv1\pmod d$ and $\a$ such that
$J_q(\a)$ has degree $\le n$ over $\Q$.  Let $H\subset
G=(\Z/d\Z)^\times\cong\gal(\Q(\mu_d)/\Q)$ be the subgroup fixing
$\Q(J_q(\a))$ so that $H$ has index $\le n$ in $G$.  Then we have
\begin{align*}
A(q,\a):=\frac{\ord_{\wp}J_q(\a)}{f}
&=\frac{1}{|H|}\sum_{t\in H}\frac{\ord_{\wp^t}(J_q(\a))}{f}\\
&=\frac{1}{|H|}\sum_{t\in H}\frac 1f \left(\sum_{i=0}^{w+1}
\sum_{j=0}^{f-1}\left\langle\frac{tp^ja_i}{d}\right\rangle-f\right)\\
&=\left(\sum_{i=0}^{w+1}\frac{1}{|HP|}
\sum_{t\in HP}\left\langle\frac{ta_i}{d}\right\rangle\right)-1
\end{align*}
where as before $P$ is the subgroup of $(\Z/d\Z)^\times$ generated by
$p$ and $HP$ is the subgroup generated by $H$ and $P$.  Reindexing
$\a$, we may assume that $d_0$ and $d_1$ are $\ge C\log d$.  Since $H$
has index $\le n$, so does $HP$ and so we get bounds on the inner sums
in the last displayed equation.  More precisely, by
Corollary~\ref{cor:alld}, the inner sum is $>\delta_n$ for
$i=2,\dots,w+1$, and by Proposition~\ref{prop:BasicEstimate}, if $d_0$
and $d_1$ are sufficiently large (so that $C\log d>C_{\epsilon,n}$),
the $i=0$ and $i=1$ terms are $>1/2-\epsilon$.  Applying this with
$\epsilon=\delta_n/4\le w\delta_n/4$, we see that for sufficiently
large $d$, $A(q,\a)\ge(w-1/2)\delta_n\ge\delta_n/2$.  This completes
the proof of part (2) of the theorem.  \qed

\section{Fermat Jacobians}

\subsection{}
Let $k$ be an arbitrary field with separable closure $\overline k$.
For each positive integer $d$ not divisible by the characteristic of
$k$ we consider the Fermat curve $F_d$ of degree $d$ over $k$ (the
zero locus of $\sum_{i=0}^2x_i^d$ in $\P^2$) and its Jacobian $J_d$.
If $A$ is an abelian variety over $k$, we say that ``$A$ appears in
$J_d$'' if there is a homomorphism of abelian varieties $A\to J_d$
with finite kernel.  We say ``$A$ appears in $J_d$ with multiplicity
$m$'' if $m$ is the largest integer such that $A^m$ appears in $J_d$.
The multiplicity with which $A$ appears in $J_d$ obviously depends
only on the $k$-isogeny class of $A$.

The following two theorems are the main results of this section.

\begin{thm}\label{thm:Jd0}
  Suppose that $k$ is a field of characteristic zero.  Then for every
  positive integer $g$, only finitely many $k$-isogeny classes of
  abelian varieties of dimension $\le g$ appear in $J_d$ as $d$ varies
  through all positive integers.  If $A$ is an abelian variety over
  $k$, then the multiplicity with which $A$ appears in $J_d$ is
  bounded by a constant depending only on the dimension of $A$.
\end{thm}

If $k$ has characteristic $p$ and $A$ is an abelian variety over $k$,
the $p$-rank of $A$ is by definition the dimension over $\Fp$ of the
group of $\overline k$-rational $p$-torsion points on $A$.  It is
known that the $p$-rank lies in the interval $[0,\dim A]$ and that it
is invariant under isogeny.

\begin{thm}\label{thm:Jdp}
  Suppose that $k$ is a field of characteristic $p>0$.  Then for every
  positive integer $g$, only finitely many $k$-isogeny classes of
  abelian varieties with positive $p$-rank and dimension $\le g$
  appear in $J_d$ as $d$ varies through all positive integers prime to
  $p$.  If $A$ is an abelian variety over $k$ with positive $p$-rank,
  then the multiplicity with which $A$ appears in $J_d$ is bounded by
  a constant depending only on the dimension of $A$.
\end{thm}

\begin{rems}
\hfill\break
\begin{enumerate}
\item \vskip-12pt We repeat that the constants in the theorems depend
  only on the dimension $g$.  In particular, they are independent of
  the characteristic of $k$.  As will be clear from the proof, they
  are also effectively computable.
\item Theorem~\ref{thm:Jd0} is already known in a more precise
  quantitative form by results of Aoki \cite{Aoki}, building on work
  of Koblitz, Rohrlich, and Shioda.  Theorem~\ref{thm:Jdp} may be
  known to experts but to my knowledge is not in the literature.  We
  will give a very simple proof of Theorem~\ref{thm:Jdp} for $k$
  finite and use this to deduce the general case and
  Theorem~\ref{thm:Jd0}.
\item It is proven in \cite{TS}, and by a different method in
  \cite{UlmerR2}, that over a field of characteristic $p$, the
  multiplicity with which a supersingular elliptic curve appears in
  $J_d$ is unbounded as $d$ varies.  Thus the last part of
  Theorem~\ref{thm:Jdp} would be false without the hypothesis of
  positive $p$-rank.  It is not clear what to expect for abelian
  varieties with $p$-rank zero which are not $\overline k$-isogenous
  to a product of supersingular elliptic curves.
\end{enumerate}
\end{rems}

The proofs of the theorems will be given in rest of this section.

\subsection{}\label{ss:oldnew}
If $d'<d$ is a divisor of $d$, then there is a canonical surjective
morphism $F_d\to F_{d'}$ ($x_i\mapsto x_i^{d/d'}$) which (because
$F_d\to F_{d'}$ is totally ramified at some place) induces an
injection of Jacobians $J_{d'}\into J_d$.  We define the {\it old
  part\/} $J_d^\text{old}$ to be the abelian subvariety of $J_d$
generated by the images of the morphisms $J_{d'}\into J_d$ as $d'$
varies through proper divisors of $d$ and we define the {\it new
  part\/} $J_d^\text{new}$ of $J_d$ to be the abelian variety over $k$
(well-defined only up to $k$-isogeny) such that $J_d$ is isogenous to
$J_d^\text{new}\times J_d^\text{old}$.  It is not hard to check, for
example by using the zeta function calculation mentioned in
\ref{ss:FermatZetas} below, that $J_d$ is isogenous to
$\prod_{d'|d}J_{d'}^\text{new}$.

Theorem~\ref{thm:Jd0} therefore follows from the statement that there
is a constant $C_g$ depending only on $g$ such that no abelian variety
$A$ of dimension $\le g$ appears in $J_d^\text{new}$ for any $d>C_g$.
Theorem~\ref{thm:Jdp} follows from the same statement with the
additional hypotheses that $A$ has positive $p$-rank and $d$ is not
divisible by $p$.

\subsection{}\label{ss:algfields}
Given a field $k$, let $\F$ be its prime field and $k_0$ be the
algebraic closure of $\F$ in $k$.  Then $k_0$ is a perfect field and
so the extension $k/k_0$ is regular.  The Fermat Jacobian $J_d$ and
its new part $J_d^\text{new}$ are defined over $\F$ and so if $A$ is
an abelian variety over $k$ which appears in
$J_d^\text{new}\times_{\F} k$ then there is an abelian variety $A_0$
defined over $k_0$ which appears in $J_d^\text{new}\times_{\F} k_0$
and with $A_0\times_{k_0} k\cong A$.  (This is an old result of Chow
which has been given a detailed modern treatment by Conrad, see
\cite{Conrad}*{3.21}.)  Moreover, the abelian variety $A_0$ and the
morphism $A_0\to J_d^\text{new}$ are both defined over some finite
extension of $\F$.  Thus it will suffice to prove the existence of the
constants $C_g$ mentioned at the end of Subsection~\ref{ss:oldnew}
(depending only on $g$, not on $k$) for the cases when $k$ is a number
field or a finite field.

\subsection{}\label{ss:HondaTate}
Let $k$ be $\Fq$, the subfield of $\Fpbar=\O_{\Qbar}/\p$ with $q$
elements.  A {\it Weil $q$-integer of weight 1\/} is an algebraic
integer $\alpha$ whose absolute value in every complex embedding is
$q^{1/2}$.  For the rest of this section we will call these simply
{\it Weil numbers\/}.

Honda-Tate theory \cite{TateHT} says that $\Fq$-isogeny classes of
$\Fq$-simple abelian varieties are in bijection with $\gal(\Qbar/\Q)$
orbits of Weil numbers.  If $A$ corresponds to $\alpha$, then
$E=\en_{\Fq}(A)\tensor\Q$ is a central simple algebra over
$\Q(\alpha)$ whose invariants in the Brauer group of $\Q(\alpha)$ can
be calculated in terms of the decomposition of $p$ in $\Q(\alpha)$.
The dimension of $A$ is $(1/2)[E:\Q(\alpha)]^{1/2}[\Q(\alpha):\Q]$ and
the eigenvalues of Frobenius on $H^1(A\times\Fqbar,\Ql)$ are the
conjugates of $\alpha$, each appearing with multiplicity
$[E:\Q(\alpha)]^{1/2}$.  The $p$-rank of $A$ is equal to the number of
eigenvalues of Frobenius which are units at $\p$ and so $A$ has
positive $p$-rank if and only if some conjugate of $\alpha$ is a unit
at $\p$.

If $C$ is a curve of genus $g$ over $\Fq$ and the $Z$-function of $C$
is
$$\frac{\prod_{i=1}^{2g}(1-\alpha_iT)}{(1-T)(1-qT)}$$ 
then the Weil numbers of the $\Fq$-simple factors of the Jacobian $J$
of $C$ are precisely the $\alpha_i$.  The multiplicity of $\alpha_i$
in the numerator is the multiplicity of the corresponding $A$ in $J$
up to $\Fq$-isogeny times $[E_A:\Q(\alpha)]^{1/2}$.

\subsection{}\label{ss:FermatZetas}
Given a positive integer $d$ and a prime power $q$ such that $q\equiv
1\pmod{d}$ we consider the Fermat curve $F_d$ over $\Fq$.  By a
theorem of Weil \cite{WeilNS}, the $Z$-function of $F_d$ over $\Fq$ is
$$\frac{\prod_{\a\in A_{d,1}}(1-J_q(\a)T)}{(1-T)(1-qT)}$$
where $A_{d,1}$ was defined in Subsection~\ref{ss:JacobiSums}.

It is clear from Weil's computation of the $Z$-function that the Weil
numbers of $J_d^\text{new}$ are precisely the $J_q(\a)$ as $\a$ runs
through
$$A'_{d,1}=\{\a=(a_0,a_1,a_2)\in A_{d,1}|\gcd(d,a_0,a_1,a_2)=1\}.$$

\subsection{The case of finite fields}
We assume $k=\Fq$ and that $A$ is an abelian variety over $k$ which
has positive $p$-rank and dimension $\le g$ and appears in
$J_d^\text{new}$.  In this case, the Weil numbers of $A$ are among the
Weil numbers of $J_d^\text{new}$.  Extending $k$ if necessary, we may
assume that $d|(q-1)$ and so the Weil numbers of $J_d^\text{new}$ are
the Jacobi sums $J_q(\a)$ where $\a\in A'_{d,1}$.  By the results
recalled in Subsections~\ref{ss:HondaTate} and \ref{ss:FermatZetas} it
follows that some $J_q(\a)$ has degree $\le 2g$ over $\Q$ and is a
unit at the prime $\p$.  This implies that $d$ is $\le C_{2g}$ where
$C_{2g}$ is the constant appearing in
Theorem~\ref{thm:GaussJacobiBounds}(2) for $n=2g$.  Therefore no
abelian variety of positive $p$-rank and dimension $\le g$ appears in
$J_d^\text{new}$ for large $d$ and this establishes
Theorem~\ref{thm:Jdp} for finite fields.  The argument in
Subsection~\ref{ss:algfields} shows that the theorem also holds for
arbitrary fields of positive characteristic.

\subsection{The case of number fields}
Suppose that $A$ is an abelian variety of dimension $\le g$ defined
over a number field $k$.  Suppose that $d$ is larger than the constant
$C(g)=C_{2g}$ appearing in Theorem~\ref{thm:Jdp} and that $A$ appears
in $J_d^\text{new}$.  Then for every prime $\wp$ of $k$ where $A$ has
good reduction, by Theorem~\ref{thm:Jdp} the reduction $A\times
\F_\wp$ has $p$-rank 0.  This would violate the following result,
which appears in \cite{Ogus}*{2.7.1}:

\begin{subsublemma} \textup{(}Katz\textup{)} 
If $A$ is an abelian variety over a number field $k$, then for
infinitely many primes of $k$, the reduction of $A$ has positive
$p$-rank.
\end{subsublemma}

For the convenience of the reader, we sketch the proof of the lemma.
Choose a prime $\ell$ larger than $2g$.  Let $L$ be a finite extension
of $k$ such that $\gal(\Qbar/L)$ acts trivially on the $\ell$-torsion
of $A$.  If $\wp$ is a prime of $L$ over the rational prime $p$ where
the reduction of $A$ has $p$-rank zero, then the trace of the
Frobenius at $\wp$ on $H^1(A\times\overline k,\Ql)$ is an integer
$\equiv0\pmod{p}$ and $\le2g(\N\wp)^{1/2}$.  If $\wp$ has absolute
degree 1 over $\Q$ (i.e., $\N\wp=p$), and $\sqrt{p}>2g$ then we see
that the trace must be zero.  On the other hand, since $\gal(\Qbar/L)$
acts trivially on $\ell$-torsion, the trace must be
$\equiv2g\pmod\ell$.  Since $\ell>2g$ this is impossible.  The
conclusion is that the reduction of $A$ at a prime of absolute degree
one over a large $p$ must have positive $p$-rank.  Such primes have
density one in $L$ and the primes under them in $k$ are an infinite
set satisfying the conclusion of the lemma.

We note that a stronger version of this result for abelian varieties
over $\Q$ is proven in \cite{BayerGonzalez}*{Prop.~5.1}.

The lemma completes the proof of Theorem~\ref{thm:Jd0} for number
fields and, as explained in Subsection~\ref{ss:algfields}, therefore
also for arbitrary fields of characteristic zero.

\section{Isotrivial abelian varieties with bounded ranks in 
$\hat\Z$ or $\Zhatp$-towers}

\subsection{}
In the rest of the paper we will give examples of abelian varieties
with bounded ranks in towers of function fields over various fields
$k$.  Before doing so, let us dispense with a trivial situation: if
$A$ is an abelian variety over $k(t)$ with good reduction away from
$0$ and $\infty$ and at worst tame ramification at $0$ and $\infty$,
then for any $d$ prime to the characteristic of $k$, the degree of the
conductor of $A$ over $k(t^{1/d})$ is bounded independently of $d$.
Geometric rank bounds then show that the rank of $A$ over $k(t^{1/d})$
is also bounded independently of $d$.  Therefore it is only
interesting to consider situations where the degree of the conductor
grows in the tower under consideration.  All our examples below are of
this type.

\subsection{}
We review some well-known facts about constant and isotrivial abelian
varieties.  Let $k$ be any field, let $L$ be the function field of a
geometrically irreducible curve $\Curve$ smooth and proper over $\spec
k$, and let $J$ be the Jacobian of $\Curve$.  Let $A_0$ be an abelian
variety over $k$ and let $A=A_0\times_kL$.  Then it is clear that
$A(L)$, the group of $L$-rational points of $A$, is canonically
isomorphic to $\mor_k(\Curve,A_0)$, the group of $k$-scheme morphisms
from $\Curve$ to $A_0$.  Moreover, we have an exact sequence
\begin{equation*}
0\to A_0(k)\to\mor_k(\Curve,A_0)\to\Hom_{k\text{-av}}(J,A_0)
\end{equation*}
where a $k$ point of $A_0$ is sent to the constant map with that value
and a morphism from $\Curve$ to $A_0$ is sent to the homomorphism of
abelian varieties induced by Albanese functoriality.  If $\Curve$ has
a $k$-rational divisor of degree 1 (for example if $k$ is finite) then
the last map above is surjective.  If $k$ is finitely generated over
its prime field, then by the Lang-N\'eron theorem, $A_0(k)$ is
finitely generated.  (See \cite{Conrad} for a modern treatment of
the Lang-N\'eron theorem.)  For any $k$, $\Hom_{k\text{-av}}(J,A_0)$ is
finitely generated and torsion free.  If $A_0$ is $k$-simple, then the
rank of $\Hom_{k\text{-av}}(J,A_0)$ is equal to the rank of the
endomorphism ring of $A_0$ times the multiplicity with which $A_0$
appears in $J$ up to $k$-isogeny.

\subsection{}
Continuing with the notation of the last subsection, suppose that
$\Curve$ is hyperelliptic, i.e., we are given a degree 2 morphism
$\Curve\to\P^1$.  Let $A'$ be the twist of $A=A_0\times_k
k(t)$ by the quadratic extension $L/k(t)$.  Since there are no
non-constant morphisms from $\P^1$ to an abelian variety, we have
$A(k(t))=A_0(k)$.  Since
\begin{equation*}
A(L)\tensor\Q\cong \left(A(k(t))\tensor\Q\right) \bigoplus 
\left(A'(k(t))\tensor\Q\right)
\end{equation*}
we conclude that $A'(k(t))$ has finite rank, bounded above by
\begin{equation}\label{eq:isotrank}
\dim_\Q A'(k(t))\tensor\Q=\dim_\Q \Hom_{k\text{-av}}(J,A_0)\tensor\Q
=\rk_\Z \Hom_{k\text{-av}}(J,A_0)
\end{equation}
with equality when $\Curve$ has a $k$-rational divisor of degree 1.

\subsection{}
We can now apply the rank formula above and our results about Fermat
Jacobians to give examples of bounded ranks in towers.  Let $K_1=k(t)$
and for every positive integer $d$ not divisible by the characteristic
of $k$, let $K_d=k(t^{1/d})$.  If the characteristic of $k$ is not 2,
let $L_1=k(u)$ with $u^2=t-1$; if the characteristic of $k$ is 2, let
$L_1=k(u)$ with $u^2+u=t$.  For all $d$ prime to the characteristic of
$k$, let $L_d=L_1K_d=k(t^{1/d},u)$.  Note that $L_d$ is the function
field of a hyperelliptic curve $C_d$ over $k$. Using ideas analogous
to \cite{UlmerR2}*{\S6}, one checks easily that there is a totally
ramified, surjective morphism from a Fermat curve $F_{n}\to\Curve_d$;
here $n=2d$ if the characteristic of $k$ is not 2 and $n=d$ if the
characteristic of $k$ is 2.  It follows that the Jacobian of $C_d$ is
an isogeny factor of $J_{n}$.  Applying the rank formula
\ref{eq:isotrank} and Theorems~\ref{thm:Jd0} and \ref{thm:Jdp} we have
the following.

\begin{thm}\label{thm:isotrivial}  
Let $k$ be a field and $A_0$ an abelian variety over $k$.  If the 
characteristic of $k$ is $p>0$, assume that $A_0$ is isogenous to a 
product of $k$-simple abelian varieties each with positive $p$-rank.  
Let $A=A\times_k k(t)$ and let $A'$ be the twist of $A$ by the quadratic 
extension $k(u)/k(t)$ where $u$ satisfies $u^2=t-1$ if the characteristic 
of $k$ is not 2 and $u^2+u=t$ if the characteristic of $k$ is 2.  Then the 
rank of the Mordell-Weil group $A'(k(t^{1/d}))$ is bounded as $d$ varies 
through all positive integers relatively prime to the characteristic of $k$.
\end{thm}

\section{Non-isotrivial elliptic curves with bounded ranks in $\Z_\ell$-towers}

\subsection{}
For examples of non-isotrivial elliptic curves with bounded ranks in
$\Z_\ell$ extensions, we consider the curve $E$ discussed in
\cite{UlmerR1} with affine equation
$$y^2+xy=x^3-t$$
over $\Fp(t)$.  

\bigskip
\begin{thm}~\label{thm:nonisol} 
Given $p$ let $S$ be the set of primes $\ell>3$ such that
$p\equiv1\pmod\ell$.  If $d$ is a product of powers of primes from
$S$, then the rank of $E(\Fpbar(t^{1/d}))$ is zero.
\end{thm}

The proof of the theorem will will be given in the rest of this section.

\subsection{}
We use the notation of Subsection~\ref{ss:JacobiSums} on Jacobi sums.
Given $p$, $d$ prime to $p$, and $\a=(a_0,\dots,a_3)\in A_{d,2}$, we
say that $\a$ is ``supersingular'' (some authors would say ``pure'')
if for one (and thus every) $q=p^f\equiv1\pmod{d}$ and all
$s\in(\Z/d\Z)^\times$ we have
$$\sum_{i=0}^3\sum_{j=0}^{f-1}\left\langle\frac{sp^ja_i}{d}\right\rangle=2f.$$
If $\a$ is supersingular, then for every prime $\wp$ of $\Q(\mu_d)$
over $p$, the valuation $\ord_{\wp}J_q(\a)$ is $f$ and this implies
that $J_q(\a)$ is a root of unity times $q$; this is the motivation
for the terminology ``supersingular''.

\subsection{}
By \cite{UlmerR1}*{6.4 and 7.7}, if $(d,6p)=1$, then the rank of
$E(\Fpbar(t^{1/d}))$ is equal to the number of elements
$t\in\Z/d\Z\setminus\{0\}$ such that $\a=(t,-6t,2t,3t)$ is
supersingular.  We are going to show that for suitable $d$ there are
no supersingular $\a$ of this form by using a descending induction
based on the following elementary identity.  Suppose that
$a\in\Z/d\Z$, $\ell$ is a prime such that $\ell^2|d$ and $\ell\nodiv
a$.  Let $H$ be the cyclic subgroup of $(\Z/d\Z)^\times$ generated by
$1+d/\ell$.  Then we have
$$\sum_{s\in H}\left\langle\frac{sa}{d}\right\rangle
=\left\langle\frac{a}{d/\ell}\right\rangle+\frac{\ell-1}{2}.$$ 

It follows that if $\a=(a_0,\dots,a_3)\in A_{d,2}$, $\ell^2|d$,
$\ell\nodiv a_i$ for all $i$, and if $\a$ is supersingular, then its
image in $A_{d/\ell}$ is also supersingular.  Indeed, we have
$$2f\ell=\sum_{s\in H}\sum_{i=0}^3
\sum_{j=0}^{f-1}\left\langle\frac{sp^ja_i}{d}\right\rangle
=\sum_{i=0}^3\sum_{j=0}^{f-1}
\left\langle\frac{p^ja_i}{d/\ell}\right\rangle+2f(\ell-1)$$
and similarly if $a$ is replaced by $ta$ with $t\in(\Z/d\Z)^\times$.

\subsection{}
We can now prove the theorem.  Suppose given $p$ and $d$ which is a
product of primes in $S$.  If the rank of $E(\Fpbar(t^{1/d}))$ were
positive, then we would have a $t\in\Z/d\Z$ such that
$\a=(t,-6t,2t,3t)$ is supersingular.  Without loss of generality we
may assume that $t\in(\Z/d\Z)^\times$ and then that $\a=(1,-6,2,3)$.
Applying the observation of the previous subsection repeatedly, we may
``reduce the level'' and find a $d'$ which is a product of distinct
primes from $S$ such that $(1,-6,2,3)\in A_{d',2}$ is supersingular.
But for such a $d'$ we have $f=1$, i.e., $p\equiv1\pmod{d'}$ and with
this one easily checks that
$$\sum_{i=0}^3\sum_{j=0}^{f-1}\left\langle\frac{p^ja'_i}{d'}\right\rangle=
\sum_{i=0}^3\left\langle\frac{a'_i}{d'}\right\rangle=1\neq
2f$$ 
and so we arrive at a contradiction to the assumption the
$E(\Fpbar(t^{1/d}))$ has positive rank.  This completes the proof of
the Theorem.

\subsection{}
The theorem shows that for any prime $p$ such that $p-1$ is not a
power of $2$ times a power of $3$, there is an elliptic curve over
$\Fp(t)$ with bounded rank in a $\Zl$ tower $\Fpbar(t^{1/\ell^n})$ for
suitable $\ell$.  We will prove a stronger result for certain small
$p$ not of this type (namely $p=2, 3, 5, 7$) in the next section and
so it seems likely that this kind of statement holds for all $p$.

In the same direction, it seems quite likely that a more refined
analysis would show that given $p$, and for $E$ as above, the rank of
$E(\Fpbar(t^{1/d}))$ is bounded as $d$ runs through all integers which 
are products of powers of primes $\ell$ such that  no power of $p$ is 
congruent to $-1$ modulo $\ell$.

Generalizing in another direction, a geometric analysis as in
\cite{UlmerR1}*{\S5} applied to the curves in \cite{UlmerR2}*{\S7}
might allow one to prove a version of Theorem~\ref{thm:nonisol} for
higher dimensional abelian varieties.

Finally, we note that it is not hard to deduce from
Theorem~\ref{thm:nonisol} that the curve defined over $\Q(t)$ by the
equation $y^2+xy=x^3-t$ has bounded rank over $\Qbar(t^{1/d})$ as $d$
ranges over all positive integers.  We omit the details since similar
results were shown by Shioda \cite{Shioda}*{Cor.~9} using closely
related techniques.

\section{Non-isotrivial elliptic curves with bounded ranks in $\Zhatp$-towers}
\subsection{}
We will use completely different techniques, unrelated to Fermat
varieties, to give a few examples of non-isotrivial elliptic curves
with bounded ranks in towers $\Fpbar(t^{1/d})$ as $d$ ranges over all
integers prime to $p$.
 
\begin{thm}\label{thm:nonisod}
If $p\in\{2,3,5,7,11\}$ then there exists an elliptic curve $E$ over
$\Fp(t)$ with $j(E)\not\in\Fp$ such that the rank of
$E(\Fpbar(t^{1/d}))$ is zero for all positive integers $d$ prime to
$p$.
\end{thm}
 
The proof of the theorem, which uses ideas from \cite{Ulmerpd}, will
be given in the rest of this section.
 
\subsection{}
Given an elliptic curve $E$ over $\Fp(t)$ with $j(E)\not\in\Fp$,
choose a non-zero invariant differential $\omega$ on $E$ and let
$\Delta=\Delta(E,\omega)$ and $A=A(E,\omega)$ be the discriminant and
Hasse invariant of $E$; the definition of the latter is reviewed in
\cite{Ulmerpd}*{\S2}.  Our assumptions imply that $\Delta$ and $A$ are
non-zero elements of $\Fp(t)$.
 
Consider the following conditions on $E$:
\begin{itemize}
\item $E$ has good or multiplicative reduction at $t=0$ and $t=\infty$.
\item At every finite non-zero place of $\Fp(t)$, $E$ obtains good
  reduction over a tamely ramified extension.
\item At every finite non-zero place $v$ of $\Fp(t)$, we have 
$$\frac{\ord_v(A)}{p-1}-\frac{\ord_v(\Delta)}{12} <\frac{1}{p}.$$
\end{itemize}
Note that the third condition is automatic at places where $E$ has
good ordinary reduction, in particular at places where $A$ and
$\Delta$ are units.  Note also that if $E$ satisfies these conditions
then it continues to satisfy them over the extensions $\Fq(t^{1/d})$
for any power $q$ of $p$ and any $d$ prime to $p$.

\subsection{}
It follows from \cite{Ulmerpd}*{Section~3 and the first sentence of
Section~6} that an elliptic curve over $\Fp(t)$ satisfying the
conditions of the previous subsection has rank 0 or 1 over any
extension $K=\Fq(t^{1/d})$. To see this, we consider the Frobenius and
Verschiebung isogenies
\begin{equation*}
\xymatrix{E\ar[r]^{Fr}&E^{(p)}\ar[r]^{V}&E}
\end{equation*}
whose composition is multiplication by $p$.  Section~3 of
\cite{Ulmerpd} computes the Selmer groups for $Fr$ and $V$ in terms of
the reduction types of $E$, $A$, and $\Delta$.  Under the conditions
of the previous subsection, the results are that $\sel(K,V)=0$ and
$\sel(K,Fr)$ is zero if $E$ has good reduction at 0 or $\infty$ and has
order $p$ if $E$ has multiplicative reduction at both 0 and
$\infty$.

We have an exact sequence
$$E^{(p)}(K)\to\sel(K,Fr)\to\sel(K,p)\to\sel(K,V)$$ 
and so the Selmer group for $p$ is either trivial or of order $p$.  In
the examples we give below, when $\sel(K,Fr)$ is non-trivial, there is
a point of order $p$ in $E^{(p)}(K)$ mapping to a generator of
$\sel(K,Fr)$ and so $\sel(K,p)=0$ and $E(K)$ has rank 0.

\subsection{}
We now give explicit examples of elliptic curves satisfying our conditions.  

Suppose $p=2$ and let $E$ be defined by
$$y^2+(t-1)xy+(t-1)^2y=x^3.$$ 
If $\omega=dx/((t-1)x+(t-1)^2)$, then $A=(t-1)$, $\Delta=t(t-1)^8$,
and $j=(t-1)^4/t$.  Standard methods show that $E$ has good, ordinary
reduction away from 0, 1, and $\infty$; that $E$ has multiplicative
reduction at 0 and $\infty$; and that at $t=1$, $E$ obtains good
reduction over an extension with ramification index 3 and the
inequality involving $A$ and $\Delta$ is satisfied.  The point
$(x,y)=((t-1)^2,(t-1)^3)$ on $E^{(2)}$ has order 2 and maps
non-trivially to $\sel(K,Fr)$ and so $\sel(K,2)=0$ for all
$K=\Fq(t^{1/d})$.

If $p=3$, let $E$ be defined by 
$$y^2=x^3+(t-1)^2x^2+t(t-1)^3x.$$
If $\omega=dx/2y$, then $A=(t-1)^2$, $\Delta=-t^2(t-1)^9$, and
$j=-(t-1)^3/t^2$.  Standard methods show that $E$ has good, ordinary
reduction away from 0, 1, and $\infty$; that $E$ has multiplicative
reduction at 0 and $\infty$; and that at $t=1$, $E$ obtains good
reduction over an extension with ramification index 4 and the
inequality involving $A$ and $\Delta$ is satisfied.  The points
$(x,y)=(t^2(t-1)^4,\pm t^2(t-1)^6)$ on $E^{(3)}$ have order 3 and map
non-trivially to $\sel(K,Fr)$ and so $\sel(K,3)=0$ for all
$K=\Fq(t^{1/d})$.

If $p=5$, let $E$ be defined by 
$$y^2=x^3+3(t-1)^4x+(t+1)(t-1)^5.$$
If $\omega=dx/2y$, then $A=(t-1)^4$, $\Delta=2t(t-1)^{10}$, and
$j=(t-1)^2/2t$.  Standard methods show that $E$ has good, ordinary
reduction away from 0, 1, and $\infty$; that $E$ has multiplicative
reduction at 0 and $\infty$; and that at $t=1$, $E$ obtains good
reduction over an extension with ramification index 6 and the
inequality involving $A$ and $\Delta$ is satisfied.  The points with
$x$ coordinate $2(t-1)^8(t^2\pm2t-1)$ on $E^{(5)}$ have order 5 and
map non-trivially to $\sel(K,Fr)$ and so $\sel(K,5)=0$ for all
$K=\Fq(t^{1/d})$.

If $p=7$, let $E$ be defined by 
$$y^2=x^3+(t-1)(t+1)^3x+5(t-1)(t+1)^5.$$
If $\omega=dx/2y$, then $A=(t-1)(t+1)^5$, $\Delta=2(t-1)^2(t+1)^9$,
and $j=4(t-1)$.  Standard methods show that $E$ has good, ordinary
reduction away from $\pm1$ and $\infty$; that $E$ has multiplicative
reduction at $\infty$; that at $t=1$, $E$ obtains good reduction over
an extension with ramification index 6 and the inequality involving
$A$ and $\Delta$ is satisfied; and that at $t=-1$, $E$ obtains good
reduction over an extension with ramification index 4 and the
inequality involving $A$ and $\Delta$ is satisfied.  It follows that
$\sel(K,7)=0$ for all $K=\Fq(t^{1/d})$.

If $p=11$, let $E$ be defined by 
$$y^2=x^3+8(t-1)(t+1)^3x+2(t-1)(t+1)^5.$$
If $\omega=dx/2y$, then $A=(t-1)^2(t+1)^8$, $\Delta=9(t-1)^2(t+1)^9$,
and $j=5(t-1)$.  Standard methods show that $E$ has good, ordinary
reduction away from $\pm1$ and $\infty$; that $E$ has multiplicative
reduction at $\infty$; that at $t=1$, $E$ obtains good reduction over
an extension with ramification index 6 and the inequality involving
$A$ and $\Delta$ is satisfied; and that at $t=-1$, $E$ obtains good
reduction over an extension with ramification index 4 and the
inequality involving $A$ and $\Delta$ is satisfied.  It follows that
$\sel(K,11)=0$ for all $K=\Fq(t^{1/d})$.

\subsection{}
The theory of modular forms modulo $p$ suggests
that the strategy employed in this section will
not work for large $p$.  Nevertheless, I conjecture that for all $p$
there are elliptic curves (indeed, absolutely simple abelian varieties
of any dimension) over $\Fp(t)$ which have bounded Mordell-Weil ranks
in the tower $\Fq(t^{1/d})$.

\end{document}